\documentclass[preprint,11pt,3p]{elsarticle}

\usepackage{amsmath,amsthm,amssymb}
\usepackage[colorlinks=true]{hyperref}
\hypersetup{urlcolor=blue, citecolor=red, linkcolor=blue}
\usepackage[capitalise,noabbrev,nameinlink]{cleveref}
\usepackage{amsfonts}
\usepackage{graphicx,epstopdf} 
\usepackage{subcaption}
\usepackage{tikz}
\usepackage{verbatim}
\usepackage{dsfont}
\usepackage{physics}
\usepackage{bm}
\usepackage[figureposition=bottom,tableposition=top]{caption}
\usepackage{multirow}

\journal{Computer Physics Communications}

\begin{document}

\begin{frontmatter}



\numberwithin{equation}{section}

\newtheorem{thm}{Theorem}

\title{ITVOLT: An Iterative Solver for the Time-Dependent Schr\"odinger Equation}


\author[inst1]{Ryan Schneider\corref{cor1}}
\ead{ryschnei@ucsd.edu}

\author[inst2]{Heman Gharibnejad}

\author[inst3]{Barry I. Schneider}

\affiliation[inst1]{
            organization={Department of Mathematics, UC San Diego},
            addressline={9500 Gilman Drive}, 
            city={La Jolla},
            state={CA},
            postcode={92093},
            country={USA}}

\affiliation[inst2]{
            organization={Computational Physics Inc.},
            addressline={8001 Braddock Road Suite 201}, 
            city={Springfield},
            state={VA},
            postcode={22151}, 
            country={USA}}
            
\affiliation[inst3]{
            organization={Applied and Computational Math Division, NIST},
            addressline={100 Bureau Drive}, 
            city={Gaithersburg},
            state={MD},
            postcode={20899}, 
            country={USA}}

\cortext[cor1]{Corresponding author}

\begin{abstract}
    We present a novel approach for solving the time-dependent Schr\"{o}dinger equation (TDSE). The method we propose converts the TDSE to an equivalent Volterra integral equation; introducing a global Lagrange interpolation of the integrand transforms the equation to a linear system, which is then solved iteratively. In this paper, we derive the method, explore its performance on several examples, and discuss the corresponding numerical details.
\end{abstract}



\begin{keyword}
Time-dependent Schr\"odinger equation \sep iterative methods \sep quadrature \sep Lagrange interpolation 



\end{keyword}

\end{frontmatter}



\vspace{.7cm}
\noindent {\bf PROGRAM SUMMARY}

\noindent
{\em Program Title:} Iterative Volterra Propagator (ITVOLT) \\
{\em CPC Library link to program files:} (to be added by Technical Editor) \\
{\em Developer's repository link:} \href{https://github.com/ry-schneider/Iterative_Volterra_Propagator.git}{https://github.com/ry-schneider/Iterative\_Volterra\_Propagator.git} \\
{\em Code Ocean capsule:} (to be added by Technical Editor)\\
{\em Licensing provisions:} MIT \\
{\em Programming language:} Modern Fortran \\
{\em Nature of problem:} ITVOLT is a solver for the time-dependent Schr\"{o}dinger equation (TDSE). More broadly, it can be applied to any problem that can be written as a Volterra integral equation. \\
{\em Solution method:} ITVOLT solves a Volterra integral equation representation of the TDSE by reducing it to a linear system via Lagrange interpolation of the integrand. It then solves the system iteratively via one of several iteration schemes. \\
{\em Additional comments including restrictions and unusual features:} General subroutines assume that the Hamiltonian for the TDSE being solved can be represented by a symmetric banded matrix whose time-dependent component is some scalar function of time (for example a pulse) multiplied by a fixed matrix. \\

\section{Introduction}
In atomic units, the time-dependent Schr\"{o}dinger equation (TDSE) is
\begin{equation}
    i \frac{\partial}{\partial t} {\bm \psi}(t) = {\bf H}(t)  {\bm\psi}(t) 
    \label{TDSE}
\end{equation}
for a Hamiltonian ${\bf H}$ and corresponding wave function ${\bm \psi}$. This first order differential equation, originally derived by Schr\"odinger \cite{TDSE}, governs the evolution of any quantum mechanical system. Solutions to the TDSE are fundamental for understanding a variety of phenomena in attosecond physics, including electron transfer \cite{GAINULLIN201772}, molecular dynamics \cite{MD}, and the interactions of particles with electromagnetic fields~\cite{TDSE_Solution_Example_1,Eberly}.\\
\indent Assuming $\bm{\psi}(t_0)$ is known, the formal solution to \eqref{TDSE} satisfies $\bm{\psi}(t) = {\bf U}(t,t_0) \bm{\psi}(t_0)$ for ${\bf U}$ a unitary, time evolution operator given by the time-ordered Dyson series
\begin{equation}
    {\bf U}(t,t_0) = 1 + \sum_{n=1}^{\infty} (-i)^n \int_{t_0}^t dt_1' \int_{t_0}^{t_1'} dt_2' \cdots \int_{t_0}^{t_{n-1}'} dt_n' {\bf H}(t_1'){\bf H}(t_2') \cdots {\bf H}(t_n').
    \label{time_evo_op}
\end{equation}
In most cases ${\bf U}$ cannot be evaluated analytically, meaning solutions to the TDSE are often only accessible numerically. In selecting an approach, we prefer a method that can solve the TDSE over large intervals at high accuracy and that is explicitly unitary, thereby conserving probability. \\
\indent A common starting point, one that is unitary, is a short-time approximation. Here, the TDSE is solved over an interval $[t_0, t_f]$ by assuming that $\hat{{\bf H}} = {\bf H}( \frac{t_0 + t_f}{2} )$, in which case \eqref{time_evo_op} simplifies and we have
\begin{equation}
   \bm{\psi}(t) = e^{-i\hat{{\bf H}}(t-t_0)} \bm{\psi}(t_0).
   \label{eqn: short time approx}
\end{equation}
To use this second-order approximation, ${\bf H}$ must be varying slowly enough on $[t_0, t_f]$ that it can be replaced by $\hat{{\bf H}}$ without losing accuracy, which in practice means the time interval must be sufficiently small. For an in-depth survey of short-time methods, each of which applies a different procedure for computing the matrix exponential, see Gharibnejad et al.\ \cite{Gharibnejad}. \\
\indent An alternative approach follows work of Magnus \cite{Magnus1954}, which showed that ${\bf U}$ can be written as the exponential of an operator $\bm{\Omega}$, itself an infinite sum of integrals of nested commutators involving ${\bf H}$. Applying this expansion requires truncating the series, evaluating the commutators and integrals, and making an approximation for the exponential of a sum of operators. Various choices for these steps produce a family of solvers of varying accuracy \cite{Trotter1959,Trotter_formula,fractal_decomp,Blanes_2009}. More general approaches for the TDSE include the $(t,t')$ method \cite{1993JChPh..99.4590P}, finite difference/finite element methods \cite{DVR}, and fourth-order Runge-Kutta (RK4) \cite{RK4TDSE}. While we are primarily interested in solvers, including ours, that can be posed in any number of spatial dimensions, we note that special emphasis has been paid to the 3D case in the literature \cite{GAINULLIN201568, pseudo}. \\
\indent Despite the many options, numerical methods for the TDSE typically struggle to achieve both accuracy and efficiency, either requiring prohibitively small time steps to obtain decent results (short-time, RK4) or making use of computationally demanding machinery regardless of step size (Magnus, $(t,t')$). The goal of this work is to develop a numerical approach that avoids these pitfalls. The method we propose begins by converting the TDSE to an equivalent Volterra integral equation. Introducing a global Lagrange interpolation of the integrand and integrating the Lagrange polynomials to obtain a set of quadrature weights reduces the problem to a set of algebraic equations, which we demonstrate can be solved iteratively.  \\
\indent The resulting method is similar in spirit to a recent one proposed by Ndong et al.\ \cite{Ndong2010ACP} and studied further by Schaefer et al.\ \cite{SCHAEFER2017368}, which also iteratively solves a Volterra integral equation representation of the TDSE. While the high level strategy of both methods is the same, our approach avoids the introduction of Chebyshev expansions and is therefore computationally simpler. \\
\indent The remainder of the paper is organized as follows. In the next section, we outline the methodology in detail. We follow that with two numerical examples, which are accompanied by practical guidelines for handling the various parameters of the method.
\section{Iterative Volterra Propagator (ITVOLT)}
\label{sec: methodology}
In this section, we develop the method and discuss associated numerical and theoretical considerations.
\subsection{Statement of the Method}
As the name suggests, our method propagates a solution to the TDSE by solving \eqref{TDSE} on intervals of the form $\left[ \tau_j, \tau_{j+1} \right]$. To do this, we start by decomposing ${\bf H}(t)$ as a sum of a time-independent operator, ${\bf H_0}$, and a time-dependent operator, ${\bf W}(t)$.  Rewriting the TDSE for $\tau_j \leq t \leq \tau_{j+1}$ as
\begin{equation}
     i \frac{\partial}{\partial t} \bm{\psi}(t) = \left[ {\bf H_0} + {\bf W} \left( \frac{\tau_j + \tau_{j+1}}{2} \right) \right] \bm{\psi}(t) + \left[ {\bf W}(t) - {\bf W} \left( \frac{\tau_j + \tau_{j+1}}{2} \right) \right] \bm{\psi}(t) ,
    \label{TDSE_rewrite}
\end{equation}
defining ${\bf H}_j = {\bf H_0} + {\bf W} \left( \frac{\tau_j + \tau_{j+1}}{2} \right)$ and ${\bf V}_j(t) = {\bf W}(t) - {\bf W} \left( \frac{\tau_j + \tau_{j+1}}{2} \right)$, and using 
\begin{equation}
    i \frac{\partial}{\partial t} \left[ e^{i {\bf H}_j t} \bm{\psi}(t) \right] = {\bf V}_j(t) \bm{\psi}(t),
    \label{tdse_intermediate}
\end{equation}
we can integrate the TDSE to obtain
\begin{equation}
    \bm{\psi}(t) =  e^{-i {\bf H}_j(t - \tau_j)} \bm{\psi}(\tau_j) - i \int_{\tau_j}^t e^{-i {\bf H}_j(t - t')} {\bf V}_j(t') \bm{\psi}(t') d t' \; \; \; \; \tau_j \leq t \leq \tau_{j+1}.
    \label{tdse_volt}
\end{equation}
This Volterra integral equation is an alternative formal solution to the TDSE. In fact, solving for $\bm{\psi}$ by repeatedly applying \eqref{tdse_volt} is equivalent to an iterative derivation of the Dyson series representation of the time-evolution operator (see \cite[\S 2]{Ndong2010ACP} for the details). \\
\indent To solve \eqref{tdse_volt} numerically, we choose a global set of points $\left\{ t_i \right\}_{i = 1}^n \in [\tau_j,\tau_{j+1}]$ and expand the integrand in the corresponding set of Lagrange polynomials $l_i(t) = \prod_{j \neq i} \frac{t - t_j}{t_i - t_j} $:
\begin{equation}
    e^{-i {\bf H}_j(t - t')} {\bf V}_j(t') \bm{\psi}(t') \approx \sum_{k=1}^n e^{-i{\bf H}_j(t - t_k)} {\bf V}_j(t_k) \bm{\psi}(t_k) l_k(t^\prime) .
    \label{lag_exp}
\end{equation}
Integrating the Lagrange polynomials over each subinterval to obtain a set of weights 
\begin{equation}
    w_{i,k} = \int_{\tau_j}^{t_i} l_k(t^\prime) dt^\prime,
    \label{quad weights}
\end{equation} 
the Volterra integral equation \eqref{tdse_volt} can then be rewritten as 
\begin{equation}
     \bm{\psi}(t_p) \approx e^{-i {\bf H}_j(t_p - \tau_j)} \bm{\psi}(\tau_j) - i \sum_{l=1}^n  w_{p,l} e^{-i {\bf H}_j(t_p - t_l)} {\bf V}_j(t_l) \bm{\psi}(t_l), 
     \label{eq:quad}
\end{equation}
at any point $t_p$. Note that the weights \eqref{quad weights} can be computed exactly using a sufficiently high-order Gauss quadrature. \\
\indent Done this way, we have a unique set of weights $\left\{ w_{i,k} \right\}_{k=1}^n$ for each quadrature point $t_i$, where $w_{i,k}$ corresponds to integrating the $k^{\text{th}}$ Lagrange polynomial over $[\tau_j, t_i]$. Moreover, the resulting quadrature is semi-global, since evaluating the integral on $[\tau_j, t_i]$ requires using all of the points, including those beyond $t_i$. \\
\indent At this point, \eqref{eq:quad} is a set of linear equations that could be solved using, for example, Gaussian elimination. For the multi-dimensional problems we are ultimately interested in solving, however, this would be computationally expensive. Instead, we focus on iterative approaches. A first option is to simply repeatedly apply \eqref{eq:quad}, where at the $(k+1)^{\text{st}}$ step we use $\bm{\psi}^{(k)}(t)$, the $k^{\text{th}}$ iterate, to obtain
\begin{equation}
    \bm{\psi}^{(k+1)}(t_p) = e^{-i {\bf H}_j(t_p - \tau_j)} \bm{\psi}^{(k)}(\tau_j) - i\sum_{l=1}^n  w_{p,l} e^{-i {\bf H}_j(t_p - t_l)} {\bf V}_j(t_l) \bm{\psi}^{(k)}(t_l) .
    \label{jacobi}
\end{equation}
For an initial approximation we can take $\bm{\psi}^{(0)}(t) = e^{-i {\bf H}_j(t - \tau_j)} \bm{\psi}(\tau_j)$, assuming $\bm{\psi}(\tau_j)$ is known. \\ 
\indent This Jacobi-type iteration uses only the values of $\bm{\psi}^{(k)}$ to evaluate $\bm{\psi}^{(k+1)}$ at any $t_p$. A more accurate but also more computationally demanding process is a Gauss-Seidel-type iteration, where at each step the currently available values of the unknowns
$\bm{\psi}^{(k+1)}(t_1), \ldots, \bm{\psi}^{(k+1)}(t_{p-1})$ are used to compute the next value according to 
\begin{equation}
    \aligned
    ({\bf I} + iw_{p,p}{\bf V}_j(t_p)) \bm{\psi}^{(k+1)}(t_p) = e^{-i {\bf H}_j (t_p - \tau_j)} &\bm{\psi}^{(k+1)}(\tau_j) - i \sum_{l=1}^{p-1} w_{p,l} e^{-i {\bf H}_j(t_p -t_l)} {\bf V}_j(t_l) \bm{\psi}^{(k+1)}(t_l) \\
    & - i \sum_{l=p+1}^n w_{p,l} e^{-i{\bf H}_j(t_p - t_l)} {\bf V}_j(t_l) \bm{\psi}^{(k)}(t_l) .
    \endaligned
    \label{tdse_gs}
\end{equation}
These iterative expressions are the basis for the first two versions of our method, which proceed as follows:
\begin{enumerate}
    \item Choose a set of quadrature points $ \left\{ t_i \right\}_{i =1}^n $ in $[\tau_j,\tau_{j+1}]$ and compute the corresponding set of weights \eqref{quad weights}. 
    \item Evaluate the inhomogeneous term $e^{-i {\bf H}_j(t - \tau_j)} \bm{\psi}(\tau_j)$ at the quadrature points, setting this to be $\bm{\psi}^{(0)}$.
    \item Given $\bm{\psi}^{(k)}$, apply either the Jacobi \eqref{jacobi} or Gauss-Seidel \eqref{tdse_gs} iteration.
    \item Continue until $\max_i ||\bm{\psi}^{(k+1)}(t_i) - \bm{\psi}^{(k)}(t_i)||_2 $ falls below a given tolerance or a maximum number of iterations is reached. 
    \item Step to the next interval $\left[ \tau_{j+1}, \tau_{j+2} \right] $, passing along the newly found value of $\bm{\psi}(\tau_{j+1})$.
\end{enumerate} 

\indent For an alternative approach, we could solve \eqref{eq:quad} via the Generalized Minimal Residual Method (GMRES) \cite{GMRES}. This popular algorithm iterates a solution to $Ax = b$ by repeatedly solving a least squares problem over Krylov subspaces of increasing dimension. Given an initial guess $x^{(0)}$ and a corresponding initial residual $r^{(0)} = b - Ax^{(0)}$, the $k$-th step of this iteration selects as $x^{(k)}$ the vector that minimizes $||Ax - b||_2$ over the Krylov subspace spanned by  $r_0, A r_0, \ldots, A^{k-1}r_0$. In practice, the least squares problem is solved by first computing an orthonormal basis for the Krylov subspace to avoid possible instability due to linear dependence. Typically, convergence is reached once the relative size of the residual $\frac{||b - Ax^{(k)}||_2}{||b||_2}$ falls below a chosen tolerance.  \\
\indent To implement GMRES here we need to rearrange \eqref{eq:quad} so that the coefficient vector $b$ corresponds to terms involving $\bm{\psi}(\tau_j)$, which we again assume is known from the previous interval. The outline for this version of our method is then the same as above with step 3 and 4 replaced by a call to an external GMRES routine. Additionally in step 2 we use the inhomogeneous term to construct the vector $b$ rather than $\bm{\psi}^{(0)}$, as GMRES typically starts with an initial guess of zero. Note that the convergence criteria for GMRES differs significantly from that of the Jacobi and Gauss-Seidel iterations; rather than stopping once successive iterates are sufficiently similar, GMRES terminates only once the corresponding residual is small enough. We comment more on the practical differences of using these iterative schemes in the next section. \\
\indent Hereafter, we refer to our method as ITVOLT, short for Iterative Volterra Propagator. When relevant, we distinguish between the three versions as ITVOLT-J (Jacobi), ITVOLT-GS (Gauss-Seidel), and ITVOLT-GMRES. \\
\indent As mentioned earlier, the high level strategy of ITVOLT was previously explored by Ndong et al.\ \cite{Ndong2010ACP}. In fact, both ITVOLT and the method they propose, referred to as Iterative Time Ordering (ITO), are special cases of the broad class of methods known as exponential integrators \cite{hochbruck_ostermann_2010}. While ITO iterates the same Volterra integral equation \eqref{tdse_volt}, its treatment of the integral is significantly more computationally demanding. In particular, Ndong et al.\ choose to expand ${\bf V}_j(t) \bm{\psi}^{(k)}(t)$ in Chebyshev polynomials, convert that expansion to a power series, and integrate the product of the power series and the matrix exponential analytically. While both Ndong et al.\ and subsequently Schaefer et al.\ \cite{SCHAEFER2017368} demonstrate that this approach is capable of high accuracy, we believe ITVOLT can do the same with a more straightforward quadrature.

\subsection{Numerical Details}
While the preceding outline emphasizes the main components of ITVOLT, there are still a number of subtle numerical details that need to be elucidated. \\
\indent First, we need to decide on a set of quadrature points. While in principle any set of points would work, the stability and accuracy of the method depends critically on how they are chosen. Equally spaced points are attractive for their simplicity, but they are subject to Runge's phenomenon \cite{runge}, which can lead to wild oscillations away from the interpolating points. Although this can be somewhat counteracted by enforcing the condition that all weights are positive \cite{HUYBRECHS2009933}, the possibility for increased instability, as we push to higher accuracy, is counterproductive. Of the many types of non-equally spaced points, we choose Gauss-Lobatto points so that the endpoints $\tau_j$ and $\tau_{j+1}$ are included. In general, increasing the number of quadrature points improves the accuracy of the quadrature but also drives up computational costs and may slow convergence.   \\
\indent Once the points are chosen, we apply a Lagrange interpolation to derive our quadrature weights. This begs the question: of all the options, why use Lagrange polynomials? The main benefit of a Lagrange interpolation is the ease with which the expansion coefficients can be found. A power expansion, for comparison, requires solving a system of equations involving the Vandermonde matrix, which is well documented to become ill-conditioned as its size grows \cite{vandermonde}. Moreover, and in contradiction to traditional lore, Lagrange interpolations exhibit stability, particularly if barycentric interpolation formulas are used \cite{barycentric, 8144747}. \\
\indent With the quadrature in place, we're left to decide between the three iterative schemes. Momentarily setting aside costs associated with applying a matrix exponential to a vector (since one iteration of Jacobi, Gauss-Seidel, or GMRES requires evaluating the same number of those), ITVOLT-J is clearly the simplest of the three; it requires additionally only matrix/vector multiplication and is parallelizable, as each $\bm{\psi}^{(k+1)}(t_i)$ can be computed independently given $\bm{\psi}^{(k)}$. Each iteration of ITVOLT-GS, on the other hand, requires solving a linear system at every quadrature point. That is, if $\bm{\psi}$ is a vector-valued function in $\mathbb{C}^d$ then \eqref{tdse_gs} is a $d \times d$ system to be solved, which may be expensive if ${\bf I} + w_{p,p} {\bf V}_j(t_p)$ is dense. ITVOLT-GMRES avoids direct system solves, but instead requires computing and storing orthonormal vectors (typically done via the Arnoldi procedure), which becomes expensive as the number of iterations grows. These observations prompt a first conclusion: ITVOLT-GS and ITVOLT-GMRES can only be more efficient than ITVOLT-J if they achieve convergence in significantly fewer iterations, thereby operating with fewer matrix exponentials overall. \\
\indent With this in mind, we might ask what kind of convergence guarantees are available for each iteration. If we use $n$ quadrature points on a $d$-dimensional problem, ITVOLT-GMRES is guaranteed to converge in exact arithmetic on any interval after at most $d(n-1)$ iterations (i.e., the size of the system being solved) \cite[Corollary 3]{GMRES}. This can be extended to finite arithmetic, where GMRES has been shown to exhibit numerical stability and superlinear convergence \cite{GMRSE_stable, VANDERVORST1993327}. \\
\indent To explore similar results for ITVOLT-J and ITVOLT-GS, consider the series of vectors
\begin{equation}
    {\bf e}_j^{(k)} = \begin{bmatrix}
        \bm{\psi}^{(k)}(t_2) - \bm{\psi}^{(k-1)}(t_2) \\
        \vdots \\
        \bm{\psi}^{(k)}(t_n) - \bm{\psi}^{(k-1)}(t_n) \\
    \end{bmatrix}, \; \; k = 1,2, \ldots 
    \label{eq: J or GS error}
\end{equation}
which record the difference in successive iterates at the quadrature points $t_2, \ldots, t_n \in [\tau_j, \tau_{j+1}]$. We omit $\bm{\psi}^{(k)}(t_1)$ from ${\bf e}_j^{(k)}$ since both the Jacobi and Gauss-Seidel iterations ensure $\bm{\psi}^{(k)}(t_1) = \bm{\psi}^{(k-1)}(t_1)$ for all $k$. Rearranging \eqref{jacobi}, we can rephrase ITVOLT-J on $[\tau_j, \tau_{j+1}]$ in terms of ${\bf e}_j^{(k)}$ as ${\bf e}_j^{(k+1)} = {\bf A}_j{\bf e}_j^{(k)}$ for 
\begin{equation}
    {\bf A}_j = -i \begin{bmatrix} 
    w_{2,2} {\bf V}_j(t_2) & w_{2,3} e^{-i {\bf H}_j(t_2-t_3)} {\bf V}_j(t_3) & \cdots & w_{2,n} e^{-i {\bf H}_j(t_2-t_n)}{\bf V}_j(t_n) \\
    w_{3,2} e^{-i {\bf H}_j(t_3 - t_2)} {\bf V}_j(t_2) & w_{3,3}{\bf V}_j(t_3) & \cdots & w_{3,n} e^{-i {\bf H}_j(t_3 - t_n)}{\bf V}_j(t_n) \\
    \vdots & \vdots & & \vdots  \\
    w_{n,2}e^{-i {\bf H}_j(t_n-t_2)} {\bf V}_j(t_2) & w_{n,3} e^{-i {\bf H}_j (t_n-t_3)} {\bf V}_j(t_3) & \cdots & w_{n,n} {\bf V}_j(t_n) \\
    \end{bmatrix}
    \label{eq: iteration matrix}
\end{equation}
Similarly, if we decompose ${\bf A}_j$ as ${\bf A}_j = {\bf L}_j + {\bf U}_j$ for ${\bf L}_j$ block lower triangular and ${\bf U}_j$ strictly block upper triangular, ITVOLT-GS can be rephrased as 
\begin{equation}
    {\bf e}_j = ({\bf I} - {\bf L}_j)^{-1} {\bf U}_j {\bf e}_j^{(k)}
    \label{eq: GS in terms of A}
\end{equation}
assuming ${\bf I} - {\bf L}_j$ is invertible. Observing $||{\bf e}_j^{(k)}||_2 \rightarrow 0$ if and only if $|| \bm{\psi}^{(k)}(t_i) - \bm{\psi}^{(k-1)}(t_i)||_2 \rightarrow 0$ for all $t_i$, we conclude that ITVOLT-J converges in exact arithmetic on $[\tau_j, \tau_{j+1}]$ if and only if $\rho({\bf A}_j)$, the spectral radius of ${\bf A}_j$, satisfies $\rho({\bf A}_j) < 1$. Meanwhile ITVOLT-GS converges as long as ${\bf I} - {\bf L}_j$ is invertible and $\rho(({\bf I} + {\bf L}_j)^{-1} {\bf U}_j) < 1$. These follow from standard results in linear algebra (see for example \cite[\S 6.10]{NLAA}). Note also that ${\bf I} - {\bf A}_j$ is the coefficient matrix used by ITVOLT-GMRES on $[\tau_j, \tau_{j+1}]$. \\
\indent In the numerical examples to come, we will explore how the spectral radius of ${\bf A}_j$ impacts the performance of both ITVOLT-J and ITVOLT-GS. We propose $\rho({\bf A}_j)$, and the many ways to cheaply estimate it, as potential tools for deciding between the iterations. Regardless, the guarantees available for Jacobi and Gauss-Seidel fall short of those for GMRES, as both ITVOLT-J and ITVOLT-GS will fail to converge in certain settings. This is part of our motivation for including ITVOLT-GMRES as an option: it serves as a fail safe, guaranteed to converge provided enough iterations are taken. Of course, a converged result is not necessarily an accurate one, and ultimately all three versions of ITVOLT are limited by the accuracy of the approximate system \eqref{eq:quad} and by extension the approximate Lagrange interpolation \eqref{lag_exp}. \\
\indent It is worth noting that GMRES comes with the option to restart after a set number of iterations as well as the option to apply various preconditioners \cite{precondition}, and these apply more broadly to ITVOLT-J and ITVOLT-GS via \eqref{eq:quad}. While they have the potential to improve convergence, we leave an exploration of preconditioning and restarting procedures to a future work and do not apply any here. \\
\indent Finally, we return to the question of evaluating matrix exponentials. In \eqref{TDSE_rewrite} above, we added and subtracted the value of ${\bf W}$ at the midpoint of $[\tau_j,\tau_{j+1}]$ in the hopes that it would better incorporate the time dependence of ${\bf H}$ into \eqref{TDSE}. By doing this, the inhomogeneous term in \eqref{tdse_volt} takes the form of a short-time approximation to $\bm{\psi}$ (although the time in the exponent may not be small). This is likely a better starting point for ITVOLT-J and ITVOLT-GS than $\bm{\psi}(\tau_j)$, but it comes at the cost of evaluating matrix exponentials, which now appear repeatedly in all three versions. In principle, any one of the many established approaches for handling matrix exponentials works with our method, and we do not specify a particular choice here. In the subsequent examples, we will discuss how ITVOLT performs with a few of the more common methods, which we describe briefly below. For a more detailed survey of these techniques, see Gharibnejad et al.\ \cite{Gharibnejad}. 
\begin{enumerate}
    \item \textbf{Full Diagonalization:} If we can diagonalize ${\bf H}_j$ as ${\bf H}_j = {\bf Q}{\bf D}{\bf Q}^T$ for ${\bf D}$ diagonal and ${\bf Q}$ orthogonal, the matrix exponentials can be computed exactly as
    \begin{equation}
        e^{-i{\bf H}_j(t-t')} = {\bf Q} e^{-i {\bf D} (t-t')} {\bf Q}^T.
    \label{diagonal_exponential}
    \end{equation}
     While this is highly accurate, it is computationally demanding and therefore impractical, even when ${\bf H}_j$ is sparse or banded.
    \item \textbf{Lanczos Iteration }\cite{Lanczos1950AnIM}\textbf{:} For an iterative approach, we can apply $e^{-i {\bf H}_j(t-t')}$ to a vector ${\bf v}$ by approximating it in a Krylov subspace. To do this, we start by building an orthonormal basis ${\bf q}_1, {\bf q}_2, \ldots, {\bf q}_m$ such that 
    \begin{equation}
        \text{span} \left\{ {\bf q}_1, {\bf q}_2, \ldots, {\bf q}_m \right\} = \text{span} \left\{ {\bf v}, {\bf H}_j {\bf v}, {\bf H}_j^2 {\bf v}, \ldots, {\bf H}_j^{m-1} {\bf v} \right\} .
        \label{lanczos_eqn1}
    \end{equation}
    Since ${\bf H}_j$ is typically symmetric, we can find the vectors ${\bf q}_i$ via a simple three-term recurrence relation, beginning with ${\bf q}_1 = {\bf v} / ||{\bf v}||_2$. Forming the matrix ${\bf Q}$ whose columns are ${\bf q}_1, \ldots, {\bf q}_m$, we then have
    \begin{equation}
        e^{-i {\bf H}_j (t-t')} {\bf v} \approx {\bf Q} e^{-i {\bf H}_j^{(m)} (t-t')} {\bf Q}^T {\bf v},
    \end{equation}
    for ${\bf H}^{(m)}$ an $m \times m$ tridiagonal matrix. Since $m$ is much smaller than the size of ${\bf H}_j$, the inner exponential can be done via diagonalization. The parameters for the method are a maximum number of iterations, a convergence criteria, and a number $l$ of vectors to re-orthogonalize against. At each step, the difference between the approximation using $k$ and $k-1$ vectors is measured; if this error does not fall below the convergence criteria another vector is added (and re-orthogonalized against the previous $l$ to improve stability), stopping once the maximum number of iterations is reached.
    \item \textbf{Chebyshev Propagation} \cite{TalEzer1984AnAA}\textbf{:} Alternatively, we can expand $e^{-i{\bf H}_j(t-t')}$ in Chebyshev polynomials $T_n$ as
    \begin{equation}
        \sum_n (2 - \delta_{0n}) e^{-i \left( \frac{\Delta}{2} + \lambda_{\min}({\bf H}_j) \right)(t-t')} J_n \left( \frac{\Delta}{2} (t-t') \right) T_n \left( -i {\bf H}_j^{\text{norm}} \right) .
        \label{cheby_eqn2}
    \end{equation}
    In this expression, $\Delta = \lambda_{\max}({\bf H}_j) - \lambda_{\min}({\bf H}_j)$ for $\lambda_{\max}({\bf H}_j)$ and $\lambda_{\min}({\bf H}_j)$ the largest and smallest eigenvalues of ${\bf H}_j$ respectively, $J_n$ is the $n^{\text{th}}$ Bessel function of the first kind, and 
    \begin{equation}
        {\bf H}_j^{\text{norm}} = 2 \frac{{\bf H}_j - \lambda_{\min}({\bf H}_j) {\bf I}}{\Delta} - {\bf I}.
        \label{eqn: normalized H}
    \end{equation}
    Note that to apply this expansion to a vector, we can exploit the recurrence relationship for the Chebyshev polynomials to save on matrix/vector multiplications. Similar to Lanczos, this method takes two parameters: a coefficient threshold and a maximum number of terms. As the expansion coefficients are computed we check their magnitude, truncating the expansion once the last coefficient falls below the threshold in magnitude or the maximum number of terms is reached.
\end{enumerate}
Since we will use them later, note that the latter two of these methods can be used in combination with a short-time approximation \eqref{eqn: short time approx} to produce two popular solvers for the TDSE, respectively Short Iterative Lanczos (SIL) and the Chebysehv Propagator. \\
\indent In the subsequent examples, we explore in detail how all of the parameters/choices available for ITVOLT affect the results of the method. At the end of the paper, we summarize our recommendations for researchers interested in using ITVOLT on a problem of their own.
\section{Numerical Examples}
In this section, we use ITVOLT to solve the TDSE for two quantum systems: a two-level atom exposed to a laser and a linearly driven harmonic oscillator. Both problems were \href{https://github.com/ry-schneider/Iterative_Volterra_Propagator}{implemented} in Fortran and compiled with gfortran version 8.5.0., and all numerical results were achieved on a Linux machine with sixteen 3.60 GHz processors. The code is original and adapted from \cite{Gharibnejad} with the exception of the GMRES routine, where we use a popular implementation of Frayss\'{e} et al.\ \cite{GMRES_CODE}. Note that, as in the previous section, all equations are presented in atomic units. \\
\indent Since both of these examples are relatively small, system solves in ITVOLT-GS (i.e., \eqref{tdse_gs}) are done directly by calling LAPACK \cite{lapack} routines for symmetric banded matrices. In principle, these could be replaced with other methods -- like GMRES or its symmetric counterpart MINRES -- to save time if the spatial dimension is large enough. Finally, we note that the implementation of ITVOLT-J is not explicitly parallel.
\subsection{Two-Level Atom}

\begin{figure}[t]
    \centering
    \includegraphics[width=.85\linewidth]{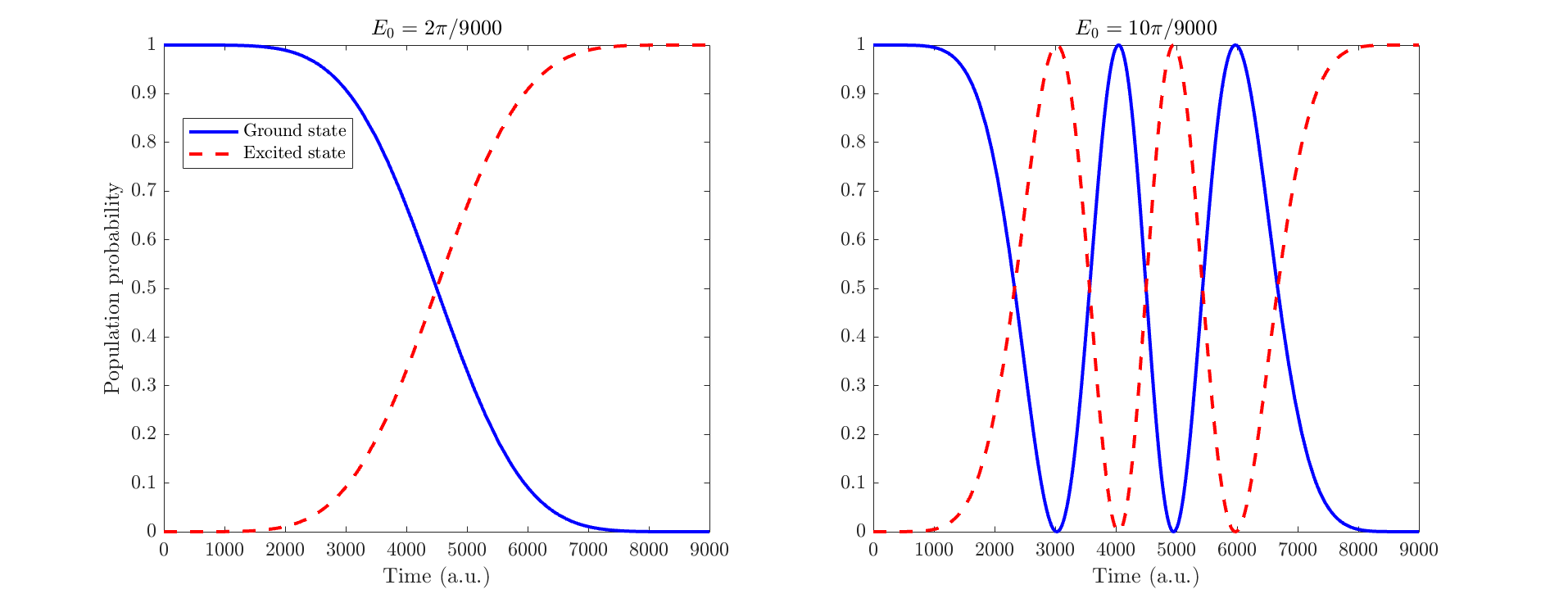}
    \caption{The ground and excited state population probabilities for the two-level atom with $T = 9000$ and $E_0 = \frac{2 \pi}{9000}$ or $E_0 = \frac{10 \pi}{9000}$, computed according to \eqref{two_level_soln}.}
    \label{two_level_data}
\end{figure}

We start with the simpler two-level problem. In the rotating wave approximation, the Hamiltonian for a two-level atom driven resonantly by a laser field $E(t)$ with transition dipole strength $1$ a.u.\ is 
\begin{equation}
    {\bf H}(t) = \begin{pmatrix} 
    0 & E(t) \\
    E(t) & 0 
    \end{pmatrix}
    \label{two_level_hamiltonian}.
\end{equation}
We solve the TDSE for this system on $[0,T]$ with a laser pulse of the form 
\begin{equation}
    E(t) =  \frac{1}{2} E_0 \sin^2 \left( \frac{ \pi t}{T} \right) 
    \label{two_level_pulse}
\end{equation}
for $E_0$ a field amplitude. Note that for this problem, since any anti-diagonal matrix can be easily diagonalized, matrix exponentials involving ${\bf H}_j = {\bf H}( \frac{\tau_j + \tau_{j+1}}{2} )$ can be done analytically. In particular,
\begin{equation}
    e^{-i {\bf H}_j (t-t')} = \begin{pmatrix} \cos(\theta) & - i \sin(\theta) \\
    - i \sin(\theta) & \cos(\theta) \\
    \end{pmatrix}, \; \; \; \theta = E \left( \frac{\tau_j + \tau_{j+1}}{2} \right)(t-t').
    \label{two_level_exp}
\end{equation} 
\indent The analytic solution corresponding to \eqref{two_level_hamiltonian} is $\bm{\psi}(t) = (c_g(t), c_e(t))^T$ for ground and excited states
\begin{equation}
    \aligned 
    & c_g(t) = \cos \left[ \frac{1}{4} E_0  \left(t - \frac{T}{2 \pi} \sin \left( \frac{2 \pi t}{T} \right) \right) \right] \\
    & c_e(t) = - i \sin \left[ \frac{1}{4} E_0 \left( t - \frac{T}{2 \pi} \sin \left( \frac{2 \pi t}{T} \right) \right) \right] .
    \endaligned 
    \label{two_level_soln}
\end{equation}
\cref{two_level_data} shows the corresponding population probabilities with $T = 9000$ and two choices of $E_0$. As anticipated, increasing the amplitude also increases the number of transitions between the two states. 

\indent To measure the accuracy of our computed solutions, we track the worst-case error in the population probabilities of each state at any propagation point $\tau_j$ used in $[0,T]$. That is, if $\bm{\psi}^{(k)}(t) = (c_1^{(k)}(t) \; \; c_2^{(k)}(t))^T $ is the converged result, we compute
\begin{equation}
    \aligned 
    & \varepsilon_{\text{sol}}^1 = \max_j \left| | c_1^{(k)}(\tau_j) |^2 - | c_g(\tau_j) |^2 \right| \; \; \text{and} \; \;  \varepsilon_{\text{sol}}^2 = \max_j \left| | c_2^{(k)}(\tau_j) |^2 - | c_e(\tau_j) |^2 \right| .
    \endaligned 
    \label{two_level_error}
\end{equation}
We also track the maximum number of iterations required to reach convergence at any step in the propagation, $k_{\max}$, as well as system run times. \\
\indent Because this problem is fairly small (and therefore quick to solve even if the number of iterations is large) we use it primarily as a comparison of the three versions of ITVOLT. With this in mind, we compute the largest spectral radius of the Jacobi iteration matrix ${\bf A}_j$ \eqref{eq: iteration matrix} over all intervals $[\tau_j, \tau_{j+1}]$:
\begin{equation}
    \rho_{\max} = \max_{j} \rho( {\bf A}_j) .
    \label{eqn: rho max}
\end{equation}
Intuitively, the value of $\rho({\bf A}_j)$ is impacted by the step size $\Delta \tau = \tau_{j+1} - \tau_j$  and the number of quadrature points $n$. We expect that $\rho({\bf A}_j)$ is large when the problem is more difficult to solve -- i.e., the step size is large or not enough quadrature points are used. \\
\indent The results for each version of ITVOLT with difference choices of $\Delta \tau$ and $n$ are shown in \cref{two_level_results} alongside the corresponding values of $\rho_{\max}$. The parameters for these runs were $E_0 = 2 \pi / 9$ and $T = 9000$, and each version of ITVOLT was allowed up to $2n - 2$ iterations, the maximum required for GMRES to reach convergence.

\renewcommand{\arraystretch}{1.4}
\begin{table}[t]
\begin{center}
\resizebox{\textwidth}{!}{
\begin{tabular}{|c|c|c||c|c|c|c|c|c|c|c|c|} 
\hline
 & & &
\multicolumn{3}{|c|}{\textbf{ITVOLT-J}}                                
&                                          
\multicolumn{3}{|c|}{\textbf{ITVOLT-GS}}
& 
\multicolumn{3}{|c|}{\textbf{ITVOLT-GMRES}} \\
\hline 
 $ \Delta \tau$ & $n$ & $\rho_{\max}$ & $\max \left\{ \varepsilon_{\text{sol}}^1, \varepsilon_{\text{sol}}^2 \right\}$ & $k_{\max}$ & $\substack{\text{System} \\ \text{time}}$ & $\max \left\{ \varepsilon_{\text{sol}}^1, \varepsilon_{\text{sol}}^2 \right\}$  & $k_{\max}$ & $\substack{\text{System} \\ \text{time}}$ & $\max \left\{ \varepsilon_{\text{sol}}^1, \varepsilon_{\text{sol}}^2 \right\}$ & $k_{\max} $ & $\substack{\text{System} \\ \text{time}}$ \\ 
 \hline
 \multirow{3}{*}{100} & 3 & $0.10$ & $5.00 \times 10^{-3}$ & 4 & 5 ms & $5.00 \times 10^{-3}$ & 3 & 4 ms & $5.00 \times 10^{-3}$ & 3 & 5 ms \\
  & 6 & $0.05$ & $1.01 \times 10^{-8}$ & 8 & 6 ms & $1.01 \times 10^{-8}$ & 5 & 9 ms & $1.01 \times 10^{-8}$ & 7 & 9 ms \\
  & 12 & $0.02$ & $1.25 \times 10^{-13}$ & 8 & 16 ms & $1.23 \times 10^{-13}$ & 5 & 19 ms & $1.87 \times 10^{-13}$ & 7 & 15 ms \\
  \hline 
  \multirow{3}{*}{500} & 6 & $1.21$ & $2.26 \times 10^{-1}$ & 10 & 7 ms & $1.71 \times 10^{-1}$ & 10 & 6 ms & $1.71 \times 10^{-1}$ & 10 & 6 ms\\
  & 12 & $0.60$ & $3.52 \times 10^{-5}$ & 22 & 15 ms & $3.52 \times 10^{-5}$ & 15 & 15 ms & $3.52 \times 10^{-5}$ & 22 & 10 ms \\
  & 24 & $0.31$ & $4.17 \times 10^{-12}$ & 25 & 41 ms & $4.54 \times 10^{-11}$ & 10 & 37 ms & $4.19 \times 10^{-13}$ & 20 & 33 ms\\
  \hline
  \multirow{3}{*}{1000} & 12 & $2.39$ & INF & - & - & INF & - & - & $5.01 \times 10^{-1}$ & 22 & 9 ms \\
  & 24 & $1.22$ & $2.72 \times 10^{-3}$ & 46 & 30 ms & $2.63 \times 10^{-3}$ & 42 & 33 ms & $2.63 \times 10^{-3}$ & 46 & 30 ms \\
  & 36 & $0.83$ & $9.80 \times 10^{-10}$ & 70 & 75 ms & $8.30 \times 10^{-10}$ & 23 & 68 ms & $8.53 \times 10^{-10}$ & 43 & 72 ms \\
  \hline 
\end{tabular}}
\end{center}
\caption{Worst-case population probability error \eqref{two_level_error} for each version of ITVOLT applied to the driven two-level atom. In all cases, a pulse with amplitude $E_0 = 2\pi / 9$ is used, the propagation is run to a final time of $T = 9000$, and $n$ quadrature points are used on intervals of size $\Delta \tau$. For each choice of $\Delta \tau$ and $n$, $\rho_{\max}$ records the maximum spectral radius of the Jacobi iteration matrix \eqref{eq: iteration matrix} over the course of the propagation. $k_{\max}$ is the maximum number of iterations required to reach convergence at any step, and listed system times are the average of five consecutive runs (rounded to the nearest millisecond). ITVOLT-J and ITVOLT-GS are run with a convergence tolerance of $10^{-10}$ while ITVOLT-GMRES is run with a tolerance of $10^{-13}$. For all three, a maximum of $2n-2$ iterations is allowed, and an error of INF indicates that the method diverged to infinity.} \label{two_level_results}
\end{table}

\indent A few observations jump out immediately from this data. First, we see consistent improvement in the results as the number of quadrature points is increased for fixed $\Delta \tau$.  When all three versions converge, there is little difference in the solution error, indicating that all are reaching the theoretical limit of the approximation being made. To reach that error, ITVOLT-GS typically takes the fewest iterations of the three, although that doesn't translate into shorter run times. This is a consequence of \eqref{two_level_exp}; since exponentials are cheap to compute, doing more of them in ITVOLT-J or ITVOLT-GMRES is preferable to the cost of solving $2 \times 2$ systems in ITVOLT-GS. \\
\indent The seventh row of the table confirms many of our observations from the previous section. When the step size is large and not enough quadrature points are used, ITVOLT-J and ITVOLT-GS both diverge to infinity while ITVOLT-GMRES converges, albeit requiring the full 22 iterations to do so. With $\rho_{\max}$ listed, we see some correlation between $\rho({\bf A}_j)$ and this behavior, suggesting the following guidelines:
\begin{enumerate}
    \item When $\rho({\bf A}_j)$ is small, all three versions are likely to converge in only a few iterations, making ITVOLT-J the most efficient.
    \item When $\rho({\bf A}_j)$ is close to one ITVOLT-GS is the best option, likely to converge in significantly fewer iterations than ITVOLT-J or ITVOLT-GMRES.
    \item For all other cases ITVOLT-GMRES is preferable, as ITVOLT-J and ITVOLT-GS may not only fail to converge but in fact diverge to infinity. 
\end{enumerate}
While in practice the spectral radius of ${\bf A}_j$ could be computed on each interval $[\tau_j, \tau_{j+1}]$ and used to choose between the iteration types, this requires forming ${\bf A}_j$ and finding its largest eigenvalue, which is computationally prohibitive. Instead, we note that $\rho({\bf A}_j)$ can be upper bounded cheaply by any consistent matrix norm of ${\bf A}_j$ or via the well-known Gershgorin circle theorem.  
\subsection{Driven Harmonic Oscillator}

\begin{figure}[t]
    \centering
    \includegraphics[width = .9\linewidth]{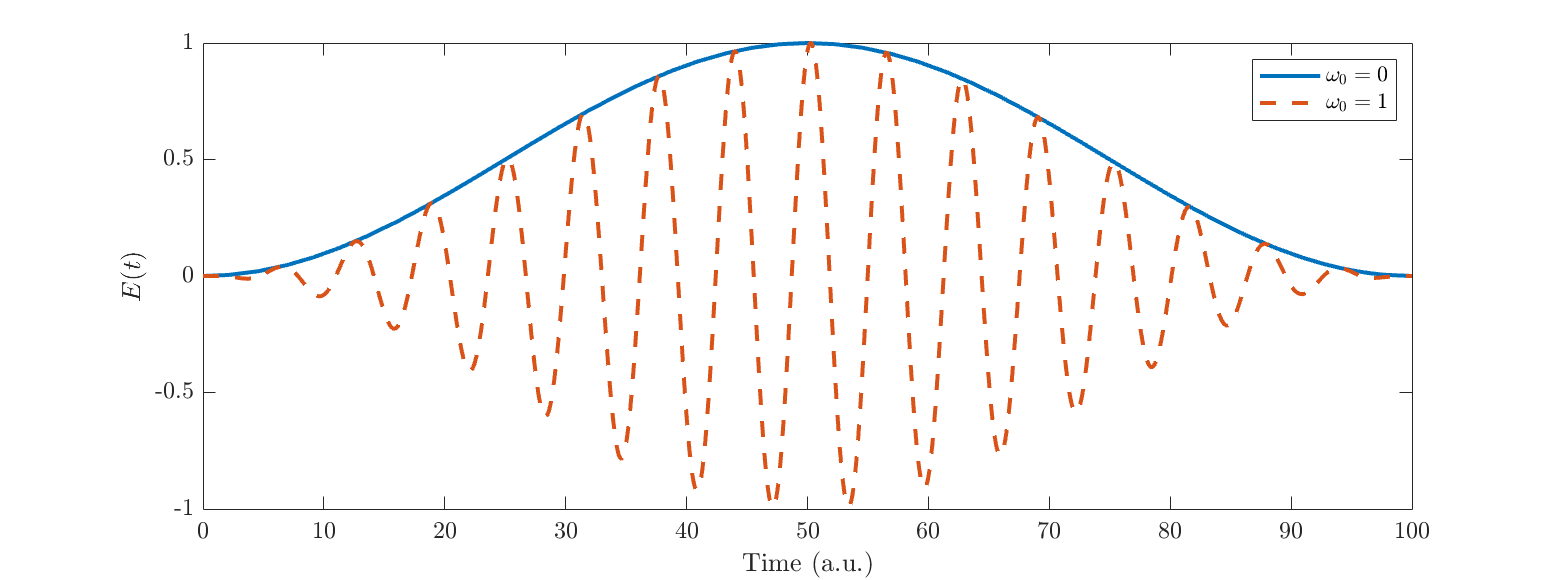}
    \caption{Driving field applied to the harmonic oscillator \eqref{oscillator_pulse} with $E_0 = 1$ and $T = 100$.}
    \label{oscillator_pulse_plot}
\end{figure}

While the results of the previous example are promising, the TDSE for a two-level atom is fairly simple. To put ITVOLT to the test on a larger, more computationally demanding problem, we next solve the TDSE for a driven harmonic oscillator with frequency $\omega = 1$ a.u. \\
\indent In atomic units, the Hamiltonian for this system is
\begin{equation}
     H(x,t) = - \frac{1}{2} \frac{ \partial^2}{\partial x^2} + \frac{1}{2} x^2 + xE(t)  
    \label{oscillator_hamiltonian}
\end{equation}
for a field
\begin{equation}
    E(t) = E_0 \sin^2 \left( \frac{ \pi t}{T} \right) \cos( \omega_0 t) .
    \label{oscillator_pulse}
\end{equation}
Here, $E_0$ is the amplitude of the pulse, $T$ is the final propagation time, and $\omega_0$ is the driving frequency. \cref{oscillator_pulse_plot} shows $E(t)$ for $E_0 = 1$, $T = 100$, and two choices of $\omega_0$. \\
\indent One of the challenges of this problem comes from the spatial variable $x$: at each time $t$, we need to decide how to represent the Hamiltonian \eqref{oscillator_hamiltonian} in space. Ndong et al.\ \cite{Ndong2010ACP}, who also treat this problem, choose to do this by representing $H$ on a Fourier grid, replacing the spatial derivative with a finite difference approximation over the grid. While our method could handle the matrix representation of $H$ obtained this way, we choose instead to expand the wave function $\psi(x,t)$ in the eigenfunctions of the unforced harmonic oscillator. That is,
\begin{equation}
    \psi(x,t) = \sum_{n = 0}^{m-1} c_n(t) \psi_n(x) 
    \label{oscillator_expansion}
\end{equation}
for $\psi_n(x)$ the $n^{\text{th}}$ eigenfunction, which solves the time-independent Schr\"odinger equation for the unforced oscillator with energy $E_n = n + \frac{1}{2}$. If ${\bf c}(t) = (c_0(t) \; \cdots \; c_{m-1}(t))^T$, this allows us to convert the TDSE to a set of coupled differential equations in time of the form
\begin{equation}
    i \frac{\partial}{\partial t} {\bf c}(t) = {\bf H}(t) {\bf c}(t),
\end{equation}
for the matrix
\begin{equation}
    {\bf H}(t) = \begin{bmatrix}
    E_0 & \frac{1}{\sqrt{2}} E(t) &  &  &\\
    \frac{1}{\sqrt{2}} E(t) & E_1 & E(t) &  & \\
    & \ddots & \ddots & \ddots & &\\ 
     & & \sqrt{\frac{m-2}{2}} E(t) & E_{m-2} & \sqrt{ \frac{m-1}{2}} E(t)  \\
    & & & \sqrt{ \frac{m-1}{2}} E(t) & E_{m-1} \\
    \end{bmatrix} 
    \label{oscillator_matrix_rep}
\end{equation}
The tridiagonal structure of ${\bf H}$ stems from $x$, which acts on the eigenfunctions $\psi_n(x)$ as a ladder operator.  \\
\indent In the case where the wave function $\psi(x,t)$ is initially in one of the eigenstates $\psi_n(t)$, the solution to the TDSE is known analytically. That is, if $\psi(0,x) = \psi_n(x)$ the TDSE has solution
\begin{equation}
    \Psi_n(x,t) = \frac{1}{\pi^{\frac{1}{4}}\sqrt{2^n n!}} e^{i p_0(t)x}e^{-i \int_0^t \delta(t') + E_n dt'}e^{-\frac{1}{2}(x - x_0(t))^2} H_n(x - x_0(t)) ,
    \label{oscillator_soln}
\end{equation}
where $x_0(t)$ and $p_0(t) = x_0'(t)$ are the position and momentum of a classical forced harmonic oscillator satisfying
\begin{equation}
    \frac{\partial^2}{\partial t^2} x_0(t) + x_0(t) = -E(t)
    \label{classical_oscillator}
\end{equation}
with initial conditions $x_0(0) = 0 = x_0'(0)$. Additionally, $\delta(t)$ is the classical Lagrangian for an unforced harmonic oscillator with the unforced position replaced by $x_0(t)$, and $H_n$ is the $n^{\text{th}}$ Hermite polynomial. Details for how this wave function is derived can be found in \cite{doi:10.1139/p58-038}. \\
\indent For our purposes we assume $\psi(0,t) = \psi_0(t)$, in which case \eqref{oscillator_soln} simplifies and the probability that the wave function is in the ground state of the unforced harmonic oscillator as a function of time can be computed exactly as 
\begin{equation}
    P_0(t) = \left| \int_{- \infty}^{\infty} \Psi_0(x,t) \psi_0(x) dx \right|^2 = \left|e^{\frac{1}{4}(x_0(t) + ip_0(t))^2 - \frac{1}{2} x_0(t)^2} \right|^2 .
    \label{oscillator_ground_prob}
\end{equation}
More generally, the probability that the system is in the $n^{\text{th}}$ state of the unforced oscillator is 
\begin{equation}
    P_n(t) = \frac{2^n}{n!} P_0(t) \left| \frac{1}{2} (x_0(t) + i p_0(t)) \right|^{2n}.
    \label{oscillator_excited_prob}
\end{equation}
To evaluate these at any time $t$ we need only find $x_0(t)$ and $p_0(t)$, which can be done analytically using the Green's function for the classical forced harmonic oscillator equation. \\
\indent Armed with these population probabilities and setting $H_0(x) = H(x,0)$, we can compute a time-dependent energy expectation value
\begin{equation}
    \langle \Psi_0(x,t) | H_0(x) | \Psi_0(x,t) \rangle = \sum_{n} P_n(t) E_n = \frac{1}{2} \left| x_0(t) + ip_0(t) \right|^2 + \frac{1}{2} . 
    \label{eqn: energy expectation}
\end{equation}
Since the energies of the unforced oscillator $E_n$ are linear in $n$, \eqref{eqn: energy expectation} gives us an idea of which states are dominating throughout the propagation. Moreover, for the probabilities \eqref{oscillator_excited_prob} it can be shown\footnote{See \ref{app: energy variance} for the details.}
\begin{equation}
    \langle \Psi_0(x,t) | H_0(x)^2| \Psi_0(x,t) \rangle - \langle \Psi_0(x,t)| H_0(x)| \Psi_0(x,t) \rangle^2 = \langle \Psi_0(x,t) | H_0(x) | \Psi_0(x,t) \rangle - \frac{1}{2},
    \label{eqn: energy variance}
\end{equation}
i.e., the variance of the energy is roughly equal to its expected value.

\begin{figure}[t]
    \centering
    \begin{subfigure}[t]{.695\linewidth}
        \centering
        \includegraphics[width= \linewidth]{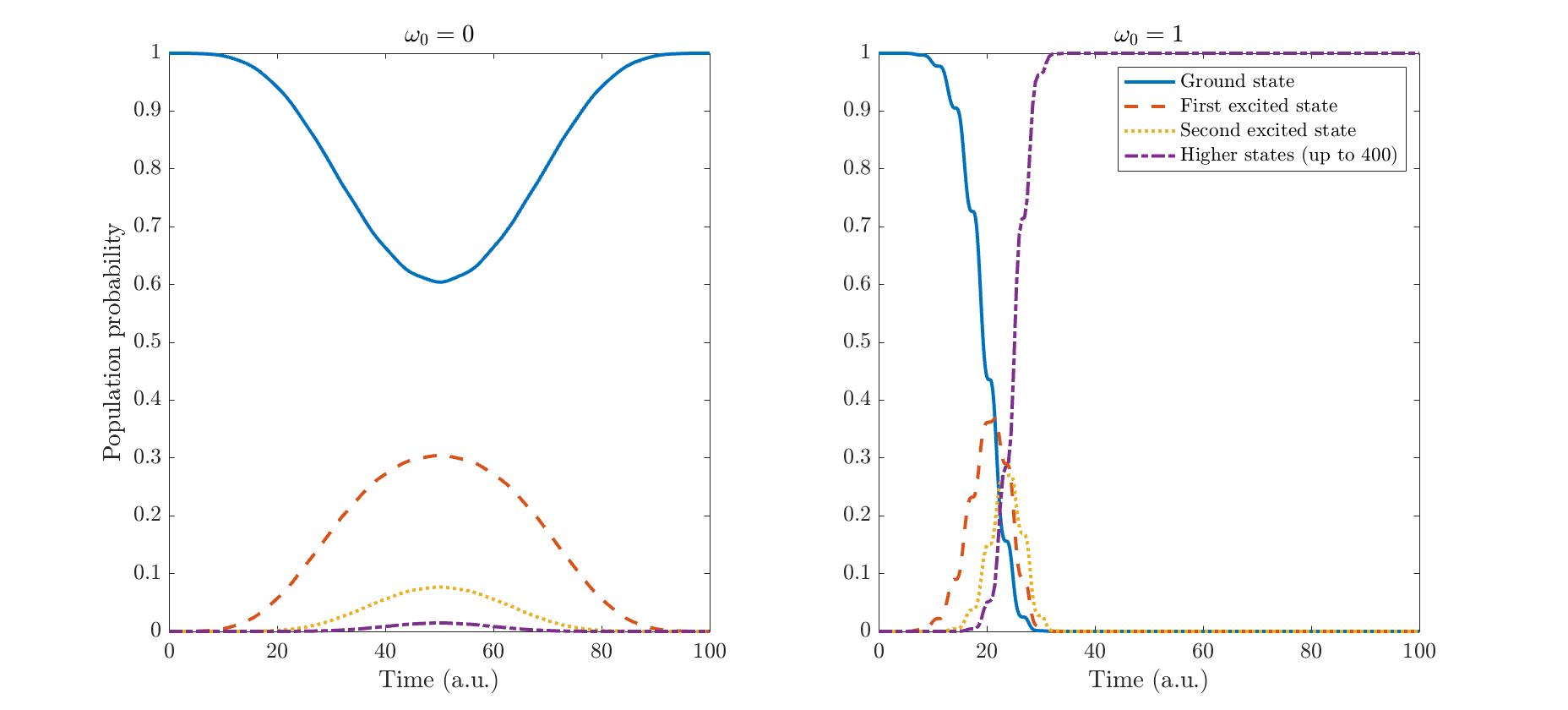}
        \caption{Population probabilities \eqref{oscillator_excited_prob} for the first 400 states.}
    \end{subfigure}%
    \begin{subfigure}[t]{.3\linewidth}
        \centering
        \includegraphics[width = \linewidth]{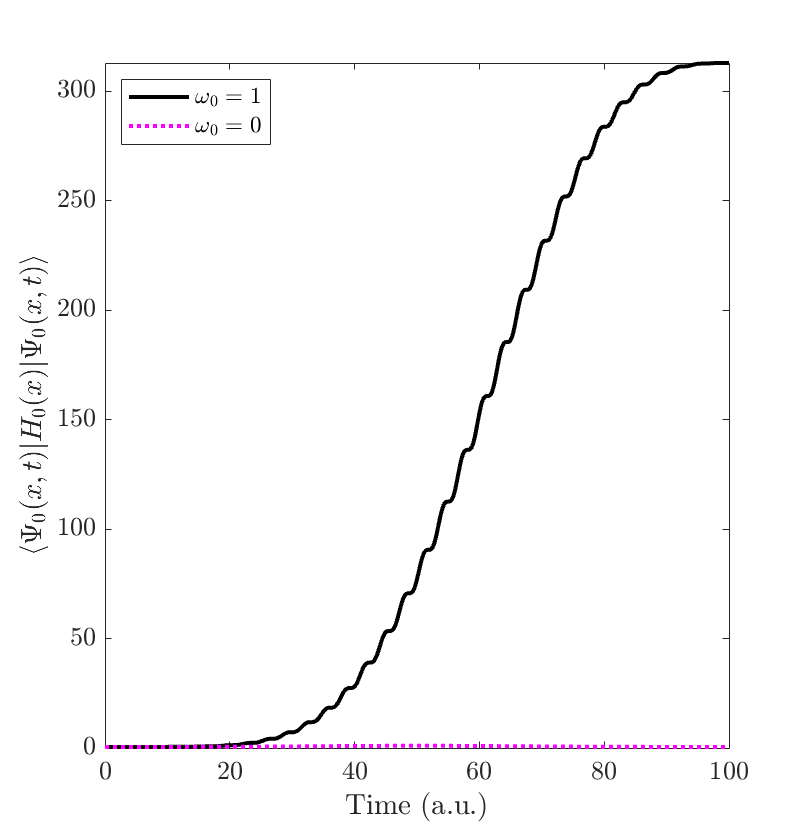}
        \caption{\centering Expected energy \eqref{eqn: energy expectation}.}
    \end{subfigure}
    \caption{Solution data for the driven harmonic oscillator with $E_0 = 1$, $T = 100$, and two choices of $\omega_0$.}
    \label{oscillator_prob_comp}
\end{figure}

\indent \cref{oscillator_prob_comp} shows the energy expectation value as well as the population probabilities for $E_0 = 1$, $T = 100$, and either $\omega_0 = 0$ or $\omega_0 = 1$. From these plots the relative difficulty of the $\omega_0 = 1$ case is clear; in this setting, the driving frequency is equal to the spacing between energy levels of the unforced oscillator, leading to a significant increase in transitions during the interaction. In light of \eqref{eqn: energy variance}, plot (b) tells us not only that the oscillator is driven to increasingly higher energies but also that the population spreads over a larger number of states as the propagation proceeds. With this in mind, we will only consider $\omega_0 = 1$ going forward. \\
\indent Similar to the previous example, we solve the TDSE on $[0,T]$ in equally spaced subintervals of the form $[\tau_j,\tau_{j+1}]$. To measure accuracy, we again compute the worst-case error in the ground state population probability at the propagation points $\tau_j$, which can be found from the converged coefficients ${\bf c}(t)$ as
\begin{equation}
    \varepsilon_{\text{sol}} = \max_j \left| |c_0(\tau_j)|^2 - P_0(\tau_j) \right|.
    \label{oscillator_error}
\end{equation}
Since this is a many state problem, we also compute the worst-case deviation from unity at any point in the propagation -- i.e.,
\begin{equation}
    \varepsilon_{\text{norm}} = \max_j \left| 1 - || {\bf c}(\tau_j)||_2^2 \right| .
    \label{oscillator_norm_error}
\end{equation}

\begin{table}[p]
\centering
\begin{subtable}{\textwidth}
\begin{center}

\begin{subtable}{\textwidth}
\begin{center}
    \resizebox{\textwidth}{!}{\begin{tabular}{|c||c|c|c||c|c|c||c|c|c|}
    \hline 
    & 
    \multicolumn{3}{|c||}{\textbf{Short Iterative Lanczos (SIL)}} &
    \multicolumn{3}{|c||}{\textbf{Chebyshev Propagator}} &
    \multicolumn{3}{|c|}{\textbf{Runge-Kutta (RK4)}} \\
    \hline
    $\Delta \tau$ & $\varepsilon_{\text{sol}}$ & $\varepsilon_{\text{norm}}$ & $\substack{\text{System} \\ \text{time}}$ & $\varepsilon_{\text{sol}}$ & $\varepsilon_{\text{norm}}$ & $\substack{\text{System} \\ \text{time}}$ & $\varepsilon_{\text{sol}}$ & $\varepsilon_{\text{norm}}$ & $\substack{\text{System} \\ \text{time}}$ \\
    \hline 
    $10^{-5}$ & $2.49 \times 10^{-10}$ & $9.68 \times 10^{-10}$ & 989.53 s & $3.86 \times 10^{-11}$ & $2.71 \times 10^{-11}$ & 13656.23 s & $3.73 \times 10^{-11}$ & $2.95 \times 10^{-11}$ & 357.79 s \\
    $10^{-4}$ & $4.09 \times 10^{-10}$ & $9.93 \times 10^{-11}$ & 123.65 s & $3.93 \times 10^{-10}$ & $1.00 \times 10^{-12}$ & 1380.32 s & $2.28 \times 10^{-12}$ & $2.93 \times 10^{-6}$ & 35.75 s \\
    $10^{-3}$ & $3.95 \times 10^{-8}$ & $1.06 \times 10^{-11}$ & 17.04 s & $3.95 \times 10^{-8}$ & $2.14 \times 10^{-12}$ & 140.19 s  & $5.45 \times 10^{-13}$ & $2.48 \times 10^{-1}$ & 3.61 s \\
    $10^{-2}$ & $3.95 \times 10^{-6}$ & $9.07 \times 10^{-13}$ & 2.96 s & $3.95 \times 10^{-6}$ & $2.36 \times 10^{-13}$ & 14.60 s & $3.27 \times 10^{-10}$ & INF & 0.36 s\\
    $10^{-1}$ & $3.95 \times 10^{-4}$ & $1.72 \times 10^{-13}$ & 1.73 s & $3.95 \times 10^{-4}$ & $2.91 \times 10^{-14}$ & 1.70 s & INF & INF & -\\
    \hline 
    
    \end{tabular}}
    \end{center}
    \caption{Error for Short Iterative Lanczos (SIL), the Chebyshev Propagator, and fourth-order Runge-Kutta (RK4). Lanczos is run with a tolerance of $10^{-12}$ for a maximum of 30 iterations, re-orthogonalizing against five vectors at each iteration. For the Chebyshev Propagator, each expansion is truncated once coefficients fall below $10^{-15}$ in modulus up to 1000 terms.}
\end{subtable}%

\vspace{5mm}

\resizebox{\textwidth}{!}{\begin{tabular}{|c|c||c|c|c||c|c|c||c|c|c|}
    \hline
    & & 
    \multicolumn{9}{|c|}{\textbf{ITVOLT with Chebyshev}} \\
    \hline
    & & 
    \multicolumn{3}{|c||}{Jacobi} & 
    \multicolumn{3}{|c||}{Gauss-Seidel} & 
    \multicolumn{3}{|c|}{GMRES} \\
    \hline 
    $\Delta \tau$ & $n$ & $\varepsilon_{\text{sol}}$ & $\varepsilon_{\text{norm}}$ & $\substack{\text{System} \\ \text{time}}$ & $\varepsilon_{\text{sol}}$ & $\varepsilon_{\text{norm}}$ & $\substack{\text{System} \\ \text{time}}$ & $\varepsilon_{\text{sol}}$ & $\varepsilon_{\text{norm}}$ & $\substack{\text{System} \\ \text{time}}$ \\
    \hline 
    \multirow{2}{*}{$10^{-2}$} & 3 & $5.82 \times 10^{-12}$ & $1.44 \times 10^{-8}$ & 29.81 s & $5.83 \times 10^{-12}$ & $1.44 \times 10^{-8}$ & 28.78 s & $6.14 \times 10^{-12}$ & $1.44 \times 10^{-8}$ & 30.10 s  \\
    & 6 & $2.84 \times 10^{-13}$ & $1.88 \times 10^{-13}$ & 108.86 s & $2.86 \times 10^{-13}$ & $1.81 \times 10^{-13}$ & 100.49 s & $2.95 \times 10^{-13}$ & $2.16 \times 10^{-12}$ & 130.95 s \\
    \hline
    \multirow{2}{*}{$10^{-1}$} & 5 & $3.29 \times 10^{-14}$ & $1.43 \times 10^{-9}$ & 19.57 s & $3.49 \times 10^{-14}$ & $1.43 \times 10^{-9}$ & 16.67 s & $1.27 \times 10^{-13}$ & $1.43 \times 10^{-9}$ & 21.93 s \\
    & 10 & $2.73 \times 10^{-14}$ & $3.60 \times 10^{-14}$ & 83.79 s & $2.52 \times 10^{-14}$ & $3.80 \times 10^{-14}$ & 69.99 s & $3.61 \times 10^{-13}$ & $1.03 \times 10^{-12}$ & 102.58 s \\
    \hline
    \multirow{2}{*}{$1.0$} & 10 & $1.30 \times 10^{-13}$ & $4.77 \times 10^{-6}$ & 82.40 s & $3.36 \times 10^{-14}$ & $4.77 \times 10^{-6}$ & 49.54 s & $1.28 \times 10^{-13}$ & $4.77 \times 10^{-6}$ & 83.19 s\\
    & 20 & $1.28 \times 10^{-13}$ & $1.01 \times 10^{-11}$ & 325.14 s & $2.48 \times 10^{-14}$ & $3.14 \times 10^{-12}$ & 188.70 s & $1.20 \times 10^{-13}$ & $2.56 \times 10^{-13}$ & 329.86 s  \\
    
    \hline 
    \end{tabular}}
    \end{center}
    \caption{Error for all three versions of ITVOLT with matrix exponentials done via the Chebyshev propagator. Chebyshev expansion parameters are the same as in (a).}
\end{subtable}%

\vspace{5mm}

\begin{subtable}{\textwidth}
\begin{center}

\resizebox{\textwidth}{!}{\begin{tabular}{|c|c||c|c|c||c|c|c||c|c|c|}
    \hline
    & & 
    \multicolumn{9}{|c|}{\textbf{ITVOLT with Lanczos}} \\
    \hline
    & & 
    \multicolumn{3}{|c||}{Jacobi} & 
    \multicolumn{3}{|c||}{Gauss-Seidel} & 
    \multicolumn{3}{|c|}{GMRES} \\
    \hline 
    $\Delta \tau$ & $n$ & $\varepsilon_{\text{sol}}$ & $\varepsilon_{\text{norm}}$ & $\substack{\text{System} \\ \text{time}}$ & $\varepsilon_{\text{sol}}$ & $\varepsilon_{\text{norm}}$ & $\substack{\text{System} \\ \text{time}}$ & $\varepsilon_{\text{sol}}$ & $\varepsilon_{\text{norm}}$ & $\substack{\text{System} \\ \text{time}}$ \\
    \hline 
    \multirow{2}{*}{$10^{-2}$} & 3 & $5.98 \times 10^{-12}$ & $1.44 \times 10^{-8}$ & 18.67 s & $5.98 \times 10^{-12}$ & $1.44 \times 10^{-8}$ & 18.46 s & $6.28 \times 10^{-12}$ & $1.44 \times 10^{-8}$ & 18.05 s \\
    & 6 & $4.36 \times 10^{-13}$ & $1.02 \times 10^{-12}$ & 125.50 s & $4.22 \times 10^{-13}$ & $1.05 \times 10^{-12}$ & 115.88 s & $4.29 \times 10^{-13}$ & $3.18 \times 10^{-12}$ & 151.76 s\\
    \hline
    \multirow{2}{*}{$10^{-1}$} & 5 & $1.45 \times 10^{-13}$ & $1.43 \times 10^{-9}$ & 54.42 s & $1.84 \times 10^{-13}$ & $1.43 \times 10^{-9}$ & 45.63 s & $1.45 \times 10^{-13}$ & $1.43 \times 10^{-9}$ & 55.63 s \\
    & 10 & $1.04 \times 10^{-13}$ & $1.49 
    \times 10^{-13}$ & 267.74 s & $3.68 \times 10^{-13}$ & $2.16 \times 10^{-13}$ & 219.79 s & $3.65 \times 10^{-13}$ & $9.73 \times 10^{-13}$ & 323.11 s  \\
    \hline 
    \end{tabular}}
    \end{center}
    \caption{Error for all three versions of ITVOLT, this time with matrix exponentials done by Lanczos. Parameters for Lanczos are again the same as in (a).}
\end{subtable}
\caption{Solution \eqref{oscillator_error} and norm \eqref{oscillator_norm_error} errors for various methods of solving the driven harmonic oscillator TDSE with $E_0 = \omega_0 = 1$, $T = 100$, and $400$ states in the expansion \eqref{oscillator_expansion}. For each ITVOLT run, $n$ quadrature points are used on intervals of size $\Delta \tau$. The convergence criteria for ITVOLT-J and ITVOLT-GS is $10^{-10}$ while ITVOLT-GMRES is run with a convergence criteria of $10^{-13}$.}
\label{table: oscillator method comparison}
\end{table}

\indent \cref{table: oscillator method comparison} summarizes these errors for all three versions of ITVOLT as well as three of the more commonly used methods for the TDSE, namely SIL, the Chebyshev Propagator, and fourth-order Runge Kutta (RK4). For a fair comparison with the two short-time methods, and to demonstrate the flexibility of ITVOLT, results for each iterative scheme are presented with exponentials done by both Lanczos and Chebyshev. Throughout, the parameters $E_0 = \omega_0 = 1$, $T = 100$, and $m = 400$ are used. Note from plot (a) of \cref{oscillator_prob_comp} that expanding in the first 400 states covers essentially the entire population over the interaction. \\
\indent For ITVOLT, many of the observations from the previous problem appear in this data. Once again, increasing the number of quadrature points used (i.e., $n$) improves the results, and the three iterations show little difference in accuracy when converged. Unlike the two-level atom, however, matrix exponentials now require calling Lanczos or Chebyshev, and we see efficiency gains for ITVOLT-GS over ITVOLT-J and ITVOLT-GMRES as a result. \\
\indent Compared with the other solution methods, ITVOLT with any iteration type and either method for handling matrix exponentials performs well. In fact, it is the only method capable of high accuracy in both measures of error simultaneously. This is the gold standard; the two short-time approximations, which we recall are unitary, demonstrate that small $\varepsilon_{\text{norm}}$ does not guarantee small $\varepsilon_{\text{sol}}$  while RK4 shows the opposite. \\
\indent The parameters used in SIL are slightly more relaxed than those used by the Chebyshev Propagator. This explains why the Chebyshev Propagator achieves better results (and also takes much longer to run) despite the fact that both are evaluating the same short-time approximation. These differences are somewhat inherited by ITVOLT, but the method is still remarkably flexible when it comes to deciding how to compute matrix exponentials. Note however that the parameters for Lanczos are too relaxed to take a step size of $\Delta \tau = 1$ in ITVOLT, hence why it is not included in the table. 

\begin{figure}[t]
    \centering
    \includegraphics[width=.9\linewidth]{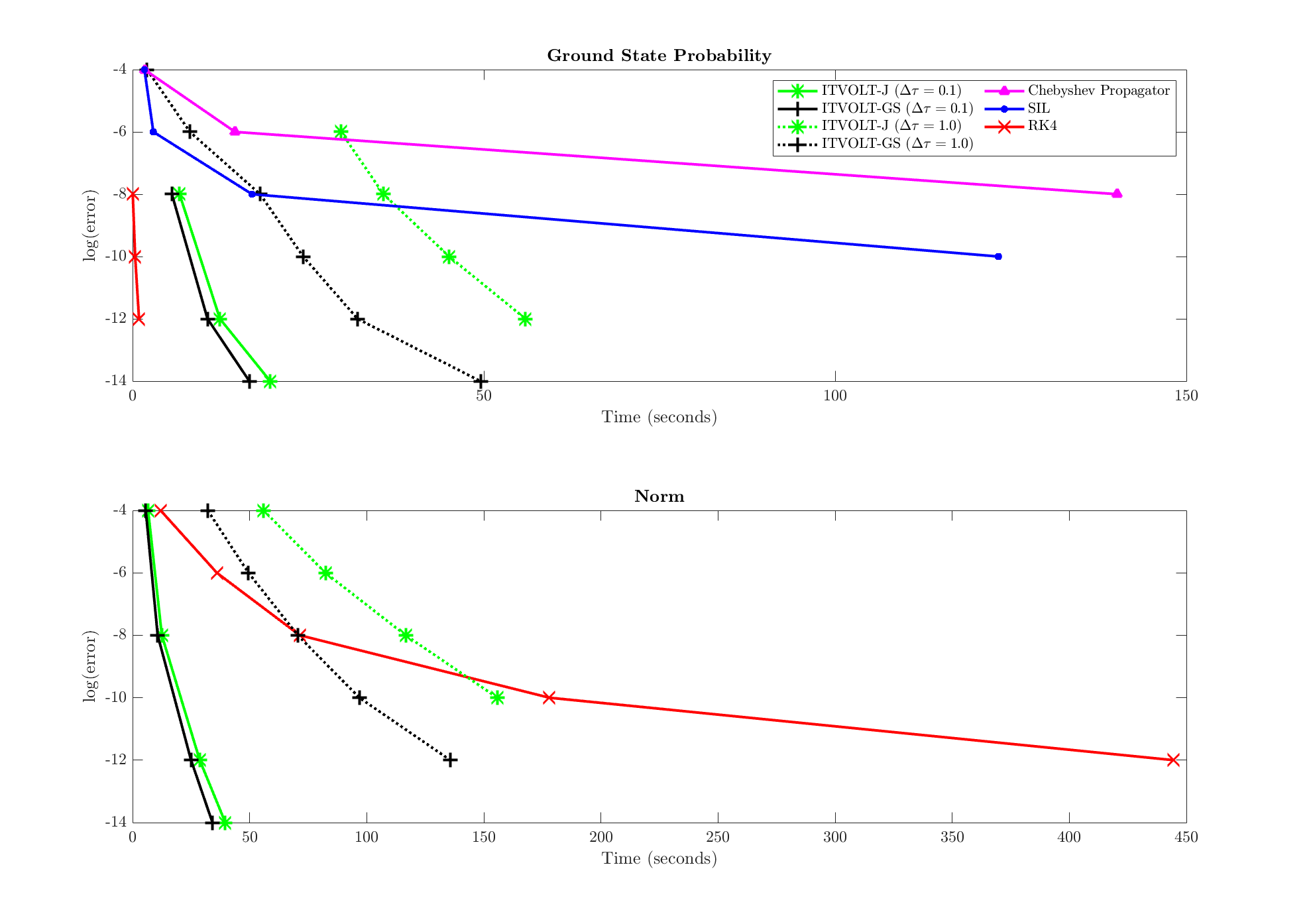}
    \caption{Minimum system time required to reach a given level of accuracy for various methods of solving the driven harmonic oscillator TDSE. Parameters are the same as in \cref{table: oscillator method comparison}, with exponentials in ITVOLT-J and ITVOLT-GS done by Chebyshev expansion. The Chebyshev Propagator and SIL are omitted from the second plot since they are explicitly unitary methods.}
    \label{fig: oscillator_method_comparison}
\end{figure}

For a more robust efficiency comparison, \cref{fig: oscillator_method_comparison} shows minimum system time 
needed by each solution method to reach a certain level of accuracy in both ground state probability $(\varepsilon_{\text{sol}})$ and unit norm $(\varepsilon_{\text{norm}})$. Errors for SIL, the Chebyshev Propagator, and RK4 are improved by shrinking the step size $\Delta \tau$. Thus, each data point for these methods corresponds to the largest value of $\Delta \tau$ that produces the desired error, where we allow only one significant digit for simplicity. Since for ITVOLT we can adjust both the step size and the number of quadrature points used, we provide two sets of points for ITVOLT-J and ITVOLT-GS, each restricting to a specific step size and adjusting the number of points used. We focus on ITVOLT-J and ITVOLT-GS with exponentials done by Chebyshev expansion as these achieve the best run times in \cref{table: oscillator method comparison}. Missing data points in either plot indicate that the method could not reach the corresponding error within the time shown or that converged results at that accuracy do not exist. \\
\indent The first plot makes clear how inefficient both SIL and the Chebyshev propagator are at computing even the ground state. This is a consequence of the short-time approximation both are based on, which is only second order. RK4 similarly struggles to achieve unit norm despite quickly and accurately computing the ground state. All three fall short of ITVOLT, which is once again the only method that can efficiently achieve high accuracy in both. Especially compelling is how quickly the error falls for both versions of ITVOLT; the difference between two consecutive data points (particularly for $\Delta \tau = 0.1$) is often only one quadrature point, which is capable of lowering the error by four orders of magnitude without significantly increasing computational costs. These plots also exhibit a nice differentiation between ITVOLT-GS and ITVOLT-J, highlighting that ITVOLT-GS not only converges faster but outperforms ITVOLT-J when $\Delta \tau$ is large.  

\begin{figure}[t]
    \centering
    \begin{subfigure}[t]{.53\linewidth}
        \centering
        \includegraphics[width= \linewidth]{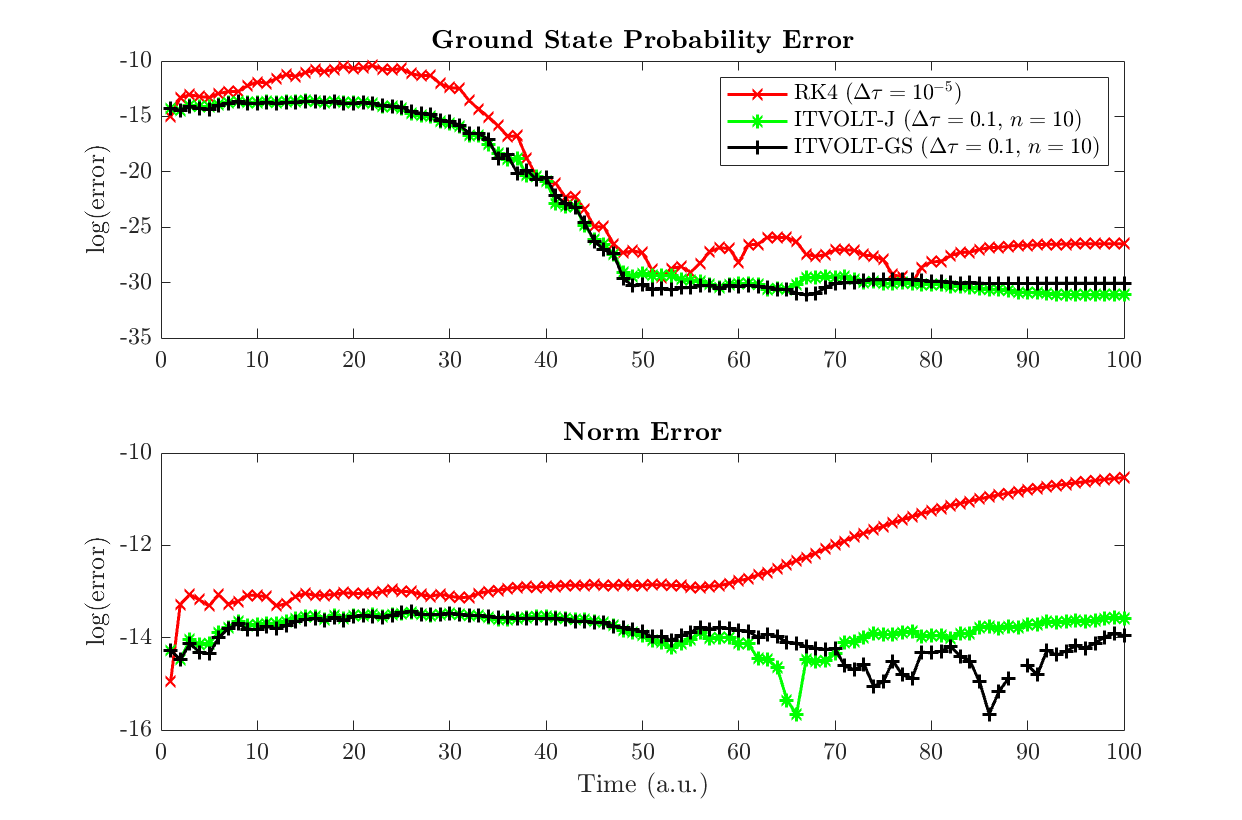}
        \caption{Evolution of error for high-accuracy methods.}
    \end{subfigure}%
    \begin{subfigure}[t]{.467\linewidth}
        \centering
        \includegraphics[width = \linewidth]{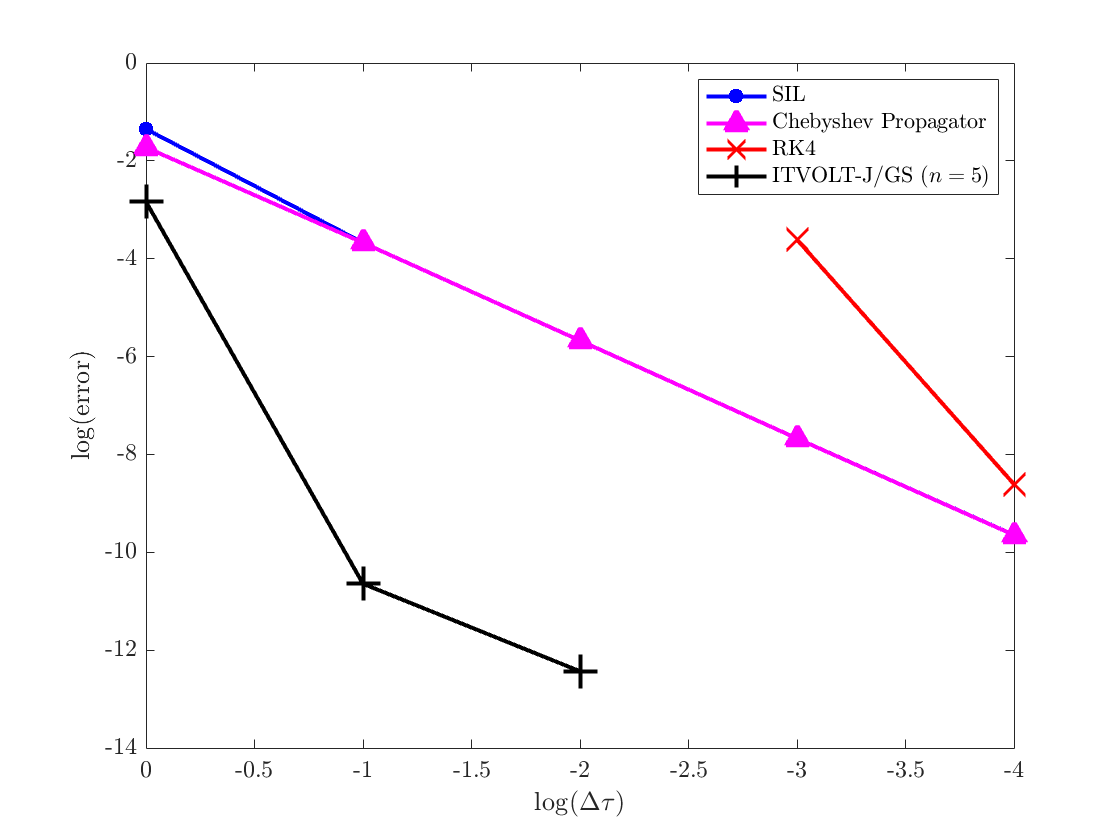}
        \caption{\centering Step size vs.\ average worst-case error in states 1-300.}
    \end{subfigure}
\caption{Solution error data for the driven harmonic oscillator TDSE with $E_0 = \omega_0 = 1$, $T = 100$, and $m = 400$. In all cases, $\Delta \tau$ is the propagation step size and $n$ is the number of quadrature points used on each interval, if relevant. Plot (a) demonstrates how error evolves over the length of the propagation for high-accuracy methods while plot (b) records average worst-case error for the first 300 states vs.\ step size for all methods. In plot (b), ITVOLT-J and ITVOLT-GS are restricted to five quadrature points and achieve the same average errors, hence are represented by one set of points. The Chebyshev Propagator and SIL are similarly represented by the same points for step sizes $\Delta \tau \leq 0.1$, while RK4 has infinite average error for $\Delta \tau \geq 0.01$. For both plots, all matrix exponentials in ITVOLT are done by Chebyshev expansion.}
    \label{fig: oscillator additional data}
\end{figure}

\indent Of course, it would be fair to ask whether focusing on worst-case errors in the ground state and norm is representative of the true performance of these methods. \cref{fig: oscillator additional data} alleviates these concerns, demonstrating that worst-case errors are not isolated and that average error over the first 300 states shows the same general trend. \\
\indent Our choice to compare ITVOLT with SIL, the Chebyshev Propagator, and RK4 is motivated by the fact that they are the most widely used methods for solving the TDSE. While ITVOLT is more complicated to implement, we believe the results presented here make a compelling case for the method. In particular, we note that because ITVOLT is built only on Lagrange interpolation, it is much less complicated to use than other methods, like ITO and $(t,t')$, that have attempted to fill the same computational gap.

\section{Conclusions and Future Work}
In this paper, we have constructed a novel iterative method for solving the time-dependent Schr\"odinger equation that is straightforward to implement and capable of accurately and efficiently solving the TDSE for a variety of systems. Based on the data presented, we make the following observations, which we hope guide researchers interested in applying ITVOLT to their own problems.
\begin{enumerate}
    \item When propagating a solution ITVOLT performs best for a medium step size, where the integral contributes significantly to the iteration but does not require too many quadrature points for accurate results.
    \item Up to a certain point, increasing the number of quadrature points used improves the results and can even reduce the number of iterations required for convergence.
    \item The spectral radius of the Jacobi iteration matrix ${\bf A}_j$ can be used to choose between the three versions of ITVOLT (as explored in the two-level problem). If all converge, ITVOLT-GS is likely to be the most efficient choice. Meanwhile, if no information is available about the difficulty of the problem ITVOLT-GMRES is the safest option.
    \item  ITVOLT is flexible when it comes to matrix exponentials, meaning a variety of methods for handling them can be used without sacrificing accuracy. 
\end{enumerate}
\indent Current efforts are focused on implementing ITVOLT on larger problems, in particular the TDSE for the three-dimensional hydrogen atom. In future work, we plan to make a more detailed comparison between ITVOLT and other high-order methods specifically for the TDSE. As mentioned earlier, we also plan to explore various preconditioned versions of the method.

\section{Acknowledgements}
The authors acknowledge support from the National Institute of Standards and Technology (NIST) in performing this research. One of the authors, R.S., is the recipient of a NIST Graduate Student Measurement Science and Engineering (GMSE) Fellowship since 2022. Additionally, a portion of this work was completed as part of a National Science Foundation Mathematical Sciences Graduate Internship (NSF MSGI). We also thank the editors and an anonymous reviewer for their helpful comments on an earlier version of the manuscript. 

\appendix
\section{Deferred Proof from Section 3.2} \label{app: energy variance}
In this appendix, we provide a proof of \eqref{eqn: energy variance}. To do this, we expand $\Psi_0(x,t)$ in the eigenfunctions $\psi_n(x)$ of $H_0(x)$ as 
\begin{equation}
    \Psi_0(x,t) = \sum_{n \geq 0} a_n(t) \psi_n(x), 
    \label{infinite expansion}
\end{equation}
where we know from \eqref{oscillator_excited_prob}
\begin{equation}
    |a_n(t)|^2 = P_n(t) = \frac{2^n}{n!}P_0(t) |h(t)|^{2n}
    \label{analytic coeff}
\end{equation}
for $h(t) = \frac{1}{2}(x_0(t) + ip_0(t))$. Recalling $E_n = n + \frac{1}{2}$, we first note
\begin{equation}
    \aligned 
    \langle \psi(x,t) | H_0(x)^2| \psi(x,t) \rangle - \langle \psi(x,t)| H_0(x)| \psi(x,t) \rangle^2 &= \sum_{n\geq 0} |a_n(t)|^2 E_n^2 - \left( \sum_{n \geq 0} |a_n(t)|^2 E_n \right)^2 \\
    & = \sum_{n \geq 0} |a_n(t)|^2n^2 - \left( \sum_{n \geq 0} |a_n(t)|^2 n \right)^2 .
    \label{eqn: variance def + energy}
    \endaligned 
\end{equation}
Now using \eqref{analytic coeff} and the fact that $\sum_{n \geq 0} |a_n(t)|^2 = 1$, we have
\begin{equation}
    \aligned 
    \sum_{n \geq 0}|a_n(t)|^2 n^2  &= \sum_{n \geq 0}|a_n(t)|^2 (n^2 - n) + \sum_{n \geq 0} |a_n(t)|^2 n \\
    & = \sum_{n \geq 2} \frac{2^n}{(n-2)!}P_0(t)|h(t)|^{2n} + \sum_{n \geq 0} |a_n(t)|^2 n \\
    & = 4|h(t)|^4 \sum_{n \geq 2} |a_{n-2}(t)|^2 + \sum_{n \geq 0}|a_n(t)|^2 n \\
    & = 4|h(t)|^4 + \sum_{n \geq 0}|a_n(t)|^2n 
    \endaligned 
\end{equation}
and similarly
\begin{equation}
    \sum_{n \geq 0} |a_n(t)|^2 n = \sum_{n \geq 1} \frac{2^n}{(n-1)!} P_0(t) |h(t)|^{2n} = 2|h(t)|^2 \sum_{n \geq 1} |c_{n-1}(t)|^2 = 2|h(t)|^2 .
\end{equation}
We conclude the proof by applying these to \eqref{eqn: variance def + energy} and noting
\begin{equation}
    \langle \Psi_0(x,t) | H_0(x) | \Psi_0(x,t) \rangle = \sum_{n \geq 0} |a_n(t)|^2 E_n = \frac{1}{2} + \sum_{n \geq 0} |a_n(t)|^2n .
\end{equation}
Since this argument depends critically on the analytic population probabilities \eqref{oscillator_excited_prob}, this kind of expectation/variance relationship is specific to our setup, where the wave function is initially in one of the eigenstates $\psi_n(x)$.

\bibliographystyle{elsarticle-num} 
\bibliography{bib.bib}

\end{document}